\documentclass[final,12pt,times,english,twoside,reqno]{amsart}

\usepackage{amsmath,amssymb,mathrsfs,amsthm,amsfonts,latexsym}
\usepackage{verbatim,setspace,enumerate}
\usepackage{matc3}
\usepackage{soul,cite,lineno,cancel,ulem,graphicx}
\usepackage{hyperref,enumitem}
\usepackage{mathtools}
\usepackage{geometry}
\newtheorem{Theorem}{Theorem}[section]
\newtheorem{Lemma}[Theorem]{Lemma}
\newtheorem{Proposition}[Theorem]{Proposition}
\newtheorem{Corollary}[Theorem]{Corollary}
\newtheorem{Remark}[Theorem]{Remark}
\newtheorem{Definition}[Theorem]{Definition}
\newtheorem{Example}[Theorem]{Example}

\newcommand{\Addresses}{{
  \bigskip
  \footnotesize \noindent
    $^{1}$   University ''Alexandru Ioan Cuza'', Faculty of Mathematics,
  Bd. Carol I, No. 11, Ia\c{s}i, 700506, ROMANIA,
  email:   croitoru@uaic.ro, OrcidID: 0000-0001-8180-3590
\\
    	$^{2}$  Petroleum-Gas University of Ploie\c{s}ti, Department of Computer Science,
      Information Technology, Mathematics and Physics,
         Bd. Bucure\c{s}ti, No. 39, Ploie\c{s}ti 100680, ROMANIA\\
          email: emilia.iosif@upg-ploiesti.ro, OrcidID: 0000-0003-0144-8811
\\
$^{3,4}$   Department of Mathematics and Computer Sciences,
     University of  Perugia -- 1, Via Vanvitelli - 06123, Perugia, ITALY,
     email: anna.sambucini@unipg.it, luca.zampogni@unipg.it ; OrcidID:0000-0003-0161-8729, 	OrcidID: 0000-0002-6047-2646
}}

\title[Some reverse inequalities for  scalar Birkhoff weak integrable ...]{
Some reverse inequalities for  scalar Birkhoff weak integrable functions}
\author{ A. Croitoru$^{1}$,  A. Iosif$^{2}$,   A. R. Sambucini$^{3}$ and L. Zampogni$^{4}$ }
\begin{document}
\begin{abstract}
The Minkowski and the H\"{o}lder inequalities play an important role in many areas of pure
and applied mathematics, such as Convex Analysis, Probabilities, Control Theory, Fixed Point
theorems, Mathematical Economics.
Also, non-additive measures, non-additive integrals and set-valued integrals are useful tools in
several areas of theoretical and applied mathematics.
In this paper we present and   prove some H\"{o}lder and Minkowski inequality (or reverse inequality) types obtained for Birkhoff weak integrable functions with respect to a non-additive measure. Then we apply these results to the interval-valued case.
\end{abstract}


\maketitle


\section{Introduction}

The classical Hölder and Minkowski inequalities are among the most powerful tools in modern mathematical analysis. 
Their influence extends from functional and convex analysis to probability theory, control theory, optimization, and mathematical economics 
(\cite{Feller, Kawabe,  radko,Milne,Yang2023,SC25}). In particular, these inequalities characterize the geometry of $L^p$-spaces 
and have numerous implications in the study of convergence, compactness, and stability of solutions to functional and integral  (\cite{CCS,Yin}).

In recent years, there has been a growing interest in generalizing such fundamental results to the setting of 
non-additive measures and non-linear integrals, motivated by the increasing role of uncertainty modeling and 
multi-criteria decision analysis. Non-additive frameworks, including fuzzy, submodular, or monotone set functions, 
offer a natural mathematical foundation for problems where the additivity assumption is no longer realistic or desirable 
(\cite{BCS,BS3,cms,torra,MS2024,Pap,GCPS}). These ideas have found applications in computer science, artificial intelligence, 
image analysis, decision theory, and statistics (\cite{mesiar1,ZZZ,W,CCGIS20}).

When $0<p<1$, the classical $L^p$-spaces cease to be normed, and the corresponding inequalities 
take a reversed or quasi-norm form. Investigating reverse Hölder and Minkowski-type inequalities in this sublinear regime 
requires nontrivial modifications of the standard additive framework. 
In particular, the Birkhoff weak integral, which is an extension of the Lebesgue integral to vector-valued or non-additive setting, introduced and studied in \cite{cdms2020,CCGIS20,cgi,CCGIS,CGI15,CG,Medjm,CGIS,cis2025},
provides a suitable tool to explore this situation.\\
Moreover
reverse inequalities are not only a  purely mathematical interest but also provide useful tools in applied mathematics, physics, engineering, and other scientific fields where classical normed-space assumptions may fail.
There are examples of their applications in analysis on metric measure spaces and geometric analysis on complex manifolds (see \cite{Frey2022});
fractional integral operators where  reverse inequalities have been applied to modified unified generalized fractional integral operators. 
 This framework is relevant in modeling viscoelastic media, anomalous diffusion, and systems with memory effects (see for example, \cite{Yang2023,Frey2022,CC,CCS,Yin,BB1,BB2}) or in 
matrix analysis and Operator Theory: here  reverse Hölder and Minkowski inequalities have been extended to matrices and quasinorms (``antinorms'' for $p<1$). These results have applications in signal processing, quantum information theory, and other areas where operators or matrices are not naturally normed. 

The aim of this paper is to establish new inequalities and reverse forms of Hölder and Minkowski inequalities 
for Birkhoff weak integrable scalar functions with respect to non-additive measures. 
These results are not only novel in the general non-additive context, but also remain new even in finitely additive settings, 
thus extending the scope of classical inequalities to broader integration theories. 
Such developments enrich the structure of generalized $L^p$-type spaces and open perspectives 
for applications in information aggregation, fuzzy control, and mathematical economics, 
where non-additive integration plays a key conceptual role.\\
The paper is organized as follows: after the Introduction, Section \ref{sec2} is devoted to  preliminaries.
Section \ref{sec3} contains some definitions,  basic results
regarding the Birkhoff weak integrability and 
 we establish some inequalities for the Birkhoff weak integral of a real function relative to a non-additive measure, such
as reverses of H\"{o}lder and Minkowski inequalities and other types of inequalities. Some applications to interval-valued or vector valued integrals  are presented in   Section \ref{sec4}. Finally a Conclusion Section follows.

 \section{Preliminaries}\label{sec2}
 Let $T$ be a non-empty set and $\mathcal{E}$ a $\sigma$-algebra of subsets of $T$.
 The integrability we consider in this paper is related to the partitions of the whole space $T$. We begin with some definitions on set functions defined on $\mathcal{E}$ and on partitions of $T$.
\begin{Definition}\label{nu}
\rm
    A set function $\nu :\mathcal{E}\rightarrow [0, \infty)$, with $\nu(\emptyset)= 0$, is called:
\begin{itemize}
\item[\ref{nu}.i)] subadditive if $\nu(A \cup B) \leq \nu(A) + \nu(B)$ , for every disjoint sets $A, B \in \mathcal{E}$.
\item[\ref{nu}.ii)] continuous  from below if for every $(B_n)_{n \in \mathbb{N}^{*}} \subset \mathcal{E}$, with $B_{n} \subset B_{n+1}$, for all $ n\in \mathbb{N}$:
\[\nu (\cup_{n=1}^{\infty} B_n) = \lim_{n \rightarrow \infty} \nu(B_n).\]
\end{itemize}
\end{Definition}
We denote by $\mathcal{M}_s$ the class of set functions $\nu:\mathcal{E}\rightarrow [0, \infty)$, with $\nu(\emptyset)= 0$, which are 
subadditive.
\begin{Example}\label{esempio1}
\rm
Let $T = \mathbb{N}$ and  $\mathcal{E}$ be its power set. 
\begin{itemize}
\item[\ref{esempio1}.i)]
Let $\nu: \mathcal{P}(\mathbb{N}) \to [0,1]$ defined by:
$\nu(A):=\sup_{n\in A}2^{-n}, \, \nu(\emptyset) =0.$
Obviously if \(A\subseteq B\), then
$\nu(A)=\sup_{n\in A}2^{-n}\le\sup_{n\in B}2^{-n}=\nu(B).$
Moreover 
for all \(A,B\subseteq\mathbb{N}\),
\[\nu(A\cup B)
=\sup_{n\in A\cup B}2^{-n}
=\max\!\left(\sup_{n\in A}2^{-n},\ \sup_{n\in B}2^{-n}\right)
=\max\{\nu(A),\nu(B)\}.
\]
Since \(\max\{x,y\}\le x+y\) for all positive \(x,y\), we get
$\nu(A\cup B)\le \nu(A)+\nu(B).$
Finally $\nu$ is continuous from below, in fact,
let \((A_k)_{k\in\mathbb{N}}\) be an increasing sequence of subsets of \(\mathbb{N}\), and let
\(A=\bigcup_{k}A_k\). Then
\[
\nu(A)
=\sup_{n\in A}2^{-n}
=\sup_{k}\sup_{n\in A_k}2^{-n}
=\sup_{k}\nu(A_k)
=\lim_{k\to\infty}\nu(A_k).
\]
The function \(\nu\) is not additive in general; for instance,
\(\nu(\{1\})=\dfrac{1}{2}\), \(\nu(\{2\})=\dfrac{1}{4}\), but
\(\nu(\{1,2\})=\dfrac{1}{2}\neq\dfrac{3}{4}\).
Thus \(\nu\) is a simple example of a monotone, subadditive, continuous from below set function.
\item[\ref{esempio1}.ii)] 
For \(A\subseteq\mathbb{N}\) define $\nu(\emptyset) = 0$  and 
\(
\nu(A):=\min \{1,\sum_{n\in A}\dfrac{1}{n}\}.
\)\\
By construction \(\nu(A)\in[0,1]\)
 for every \(A\subseteq\mathbb{N}\), and \(\nu(\emptyset)=0\).
If \(A\subseteq B\) then 
\(
\nu(A)=\min  \{1,\sum_{n\in A}\dfrac{1}{n}\}
\le
\min  \{1,\sum_{n\in B}\dfrac{1}{n}\}=\nu(B).
\)
For any \(A,B\subseteq\mathbb{N}\) 
set \(s_A:=\sum_{n\in A}\dfrac{1}{n}\) and 
\(s_B:=\sum_{n\in B}\dfrac{1}{n}\).
Then \(s_{A\cup B}\le s_A + s_B\), and 
\(
\min \{1,s_{A\cup B}\}\le\min \{1,s_A+s_B\}
\)
holds.
Therefore
\begin{eqnarray*}
\nu(A\cup B) &=& \min \{1,\sum_{n\in A\cup B}\dfrac{1}{n}\}
\le
\min \{1,\sum_{n\in A}\dfrac{1}{n}+\sum_{n\in B}\dfrac{1}{n}\}
\le 
\nu(A)+\nu(B),
\end{eqnarray*}
so \(\nu\) is subadditive.
Finally
let \((A_k)_{k\in\mathbb{N}}\) be an increasing sequence of subsets of \(\mathbb{N}\) (i.e., \ for every $k \in \mathbb{N}$, \ \(A_k\subseteq A_{k+1}\)) and set \(A=\bigcup_{k}A_k\).
Then the sequence \( \big(\sum_{n\in A_k}\dfrac{1}{n}\big)_k \) is nondecreasing  with limit \(s_A=\sum_{n\in A}\dfrac{1}{n}\) (possibly \(+\infty\)). Since the map \(x\mapsto\min \{1,x\}\) is continuous on \([0,\infty]\), we get

\[\nu(A)=\min \{1,s_A\}=\lim_{k\to\infty}\min \{1,s_{A_k} \}=\lim_{k\to\infty}\nu(A_k),\]

so \(\nu\) is continuous from below.
\end{itemize}
\end{Example}

\begin{Definition}
\rm
A property $(P)$ holds {\it $\nu$-almost everywhere} (denoted by
$\nu$-a.e.) if there exists $B \in\mathcal{E}$, with $\nu(B) = 0$, so that the property $(P)$ is valid on $T\setminus B.$
\end{Definition}

\begin{Definition}\label{uno}
\rm
Suppose $\mbox{card}(T)\geq \aleph_{0}$ (where $\mbox{card}(T)$ is the cardinality of $T$).
\begin{itemize}
\item[\ref{uno}.i)]  A countable family of nonvoid sets $P=\{B_{n}\}_{n\in \mathbb{N}}$ $\subset \mathcal{E}$
such that $\bigcup\limits_{n\in \mathbb{N}}B_{n}=T$ with $B_{i}\cap B_{j}=\emptyset ,$ when $i\neq j, \,\,  i, j \in \mathbb{N}$,
is called a (measurable) countable partition of $T$.\\
Denote by $\mathcal{C}$ the set of all countable partitions of $T$ and by $\mathcal{C}_{B}$
the set of countable partitions of $B\in \mathcal{E}$.\\

\item[\ref{uno}.ii)] For every $P$ and $P^{\prime }\in \mathcal{C}$,
$P^{\prime }$ is called finer than $P$
(denoted by $P^{\prime }\geqslant P$\,  or \, $P\leqslant P^{\prime }$) if every set of $
P^{\prime }$ is included in some set of $P$.\\

\item[\ref{uno}.iii)]   For every $P$ and $P^{\prime }\in \mathcal{C}$,\,
$P=\{B_{n}\}, P^{\prime}=\{C_{m}\}$, \, the common refinement of
$P$ and $P^{\prime }$ is defined to be the countable partition $\{B_{n}\cap C_{m}\}$,
denoted by $P\wedge P^{\prime }$.
\end{itemize}
\end{Definition}
\noindent
\section{Birkhoff weak integrability and related inequalities}\label{sec3}
This section contains definitions and basic results on the Birkhoff weak integrability and new results on reverse inequalities.
In the sequel  $T$ is a non-empty set, with $\mbox{card}(T)\geq \aleph_{0}$,\, $\mathcal{E}$ is a $\sigma$-algebra
of subsets of $T$
and $\nu:\mathcal{E}\rightarrow [0, \infty)$ is a non-negative set function, such that $\nu(\emptyset)= 0.$\\

We recall the following definition:
\begin{Definition}
\rm 
(\cite{cgi}) \label{def-int}
A  function $u:T\rightarrow \mathbb{R}$ is said to be Birkhoff weakly
 integrable on $T$  with respect to $\nu$ ($B_w-\nu$-integrable),
 if  $b\in \mathbb{R}$ exists
such that for every $\varepsilon>0$, there are
$P_{\varepsilon}\in \mathcal{C}$ and $n_{\varepsilon}\in \mathbb{N}$  
such that for every $  P \in \mathcal{C}$, $P=(B_{n})_{n\in \mathbb{N}}$, $P\geqslant P_{\varepsilon}$
and every $t_{n}\in B_{n},\,  n\in \mathbb{N}$:
\[ \Big\vert \sum_{k=1}^{n}u(t_{k})\nu(B_{k})-b \Big\vert<\varepsilon, \quad \mbox{ for every } n\geq n_{\varepsilon}.\]
$b$ is denoted by $(B_w)\int_{T}ud\nu$ or simply $\int_{T}ud\nu$ and is called
the Birkhoff weak integral of $u$ on $T$ with respect to $\nu$.
\end{Definition}

We denote by $B_w(\nu,T)$ the family of all $B_w-\nu$-integrable functions on $T$.
The Birkhoff weak integrability on every set $E \in \mathcal{E}$ is defined in the usual way.\\

In particular, by \cite[Theorem 4.2]{CGI15} $u$ is $B_w-\nu$-integrable  on $E \in \mathcal{E}$ if and only if 
$u \cdot 1_E \in B_w(\nu,T)$ and
 $(B_w)\int_{E}ud\nu =  (B_w)\int_{T}u \cdot 1_E d\nu$, where $1_E$ is the characteristic function of $E$.\\
 With the symbol $B_w(\nu)$ we denote the family of scalar functions that are  $B_w-\nu$-integrable  on every $E \in \mathcal{E}$.
  The family of $B_w-\nu$-integrable functions is closed with respect to the order of $\mathbb{R}$, in fact:
 \begin{Theorem}\label{Tix}\mbox{\rm (\cite[Corollary 3.3]{CG})}
Let $u, v:T\rightarrow \mathbb{R}$, such that $u, v\in B_w(\nu)$. 
Then $\min\{u, v\}$ and $\max\{u, v\}$ are in $B_w(\nu).$
\end{Theorem}
For other results on this topic we refer to \cite{CG,cgi,CCGIS,CGIS}.\\

In general for the gauge integrals, like Henstock, McShane, Birkhoff integrals,  no measurability condition is asked a priori, see for example \cite{BS, CS, CCGS,CCGIS, GIC, PS,KHB,SC}. \\

For the study of the  inequalities object of this research, we need sometimes the measurability of functions $u:T\rightarrow \mathbb{R}$, when it will be necessary we specify it. \\
\begin{itemize}
\item
 We denote by $\mathscr{F}(T, \mathbb{R})$ the space of all  measurable functions from $T$ to $\mathbb{R}$.\\
\end{itemize}
If $p\in(0, \infty)$ and $u:T\rightarrow \mathbb{R}$ is a function with $\vert u\vert ^{p}\in B_w(\nu, T)$, we denote, as usual,
\begin{equation}\label{*}
 \Vert u\Vert_{p}=\left( (B_w)\int_{T}\vert u\vert ^{p} d\nu \right)^{\frac{1}{p}}.
\end{equation}

Recall that two indices \(p,q>0\) are said to be conjugate if 
\mbox{ \(  p^{-1}+q^{-1}= 1.\)}
This requires \(p\geq 1\), which forces \(q\geq 1\) as well.  

Here, however, we allow \(0<p<1\). In this case the condition 
\(  p^{-1}+q^{-1}= 1\)
implies that \(q<0\).  
We will  refer to such pairs \(p,q\) again as conjugate indices and we
 will show later that the function \(\Vert \cdot \Vert : B_w(\nu, T) \to \mathbb{R}\) satisfies the 
triangle inequality when \(p \geq 1\), under suitable conditions, whereas it fails to satisfy it when 
\(0<p<1\).\\
 Consequently, in this latter case it is not a norm, and the 
vector space \(B_w^{p}(\nu, T)\) is not normed. So, in order to approach 
 the inequalities of  Minkowski and  Hölder and their reverse inequalities we consider:

\begin{Definition}\mbox{\rm  (\cite[Definition 3.5]{CG}) }\label{def-nu-int}
A set function $\nu:\mathcal{E}\rightarrow [0,\infty)$ is called {\it $\mathcal{E}$-integrable }
  if for all $B\in \mathcal{E}$, \, $1_{B}\in B_w(\nu, T)$.
\end{Definition}
Obviously, any measure $\nu:\mathcal{E}\rightarrow [0,\infty)$ is $\mathcal{E}$-integrable, see for example  (\cite{C,Medjm,CG}).
Moreover,  if we consider the set functions  given in Example \ref{esempio1}, we have
\begin{Example}\label{esempio2} \rm 
As in Example \ref{esempio1}, let $T = \mathbb{N}$ and  $\mathcal{E}$ be its power set.
 \begin{itemize}
 \item[\ref{esempio2}.i)] Let $\nu$ as in Example \ref{esempio1}.i) 
 and  take $u : \mathbb{N} \to \mathbb{R}$. 
Since, by \cite[Example II]{cgi},  
\[(B_w)\int_{T}ud\nu = \sum_{n=1}^{\infty} u(n) \nu(\{n\})\]
if the involved  series is absolutely convergent, we have that, for every $A \subset \mathbb{N}$,
\begin{eqnarray*}
\nu(A) \leq \sum_{n \in A} \nu(\{n\})  \quad \mbox{and} \quad 
(B_w)\int_{T} 1_A d\nu =\sum_{n \in A} \nu(\{n\}) = \sum_{n \in A} \big( \frac12 \big)^n <+\infty.
\end{eqnarray*}
So $\nu$ is integrable  in the sense of Definition \ref{def-nu-int}.
\item[\ref{esempio2}.ii)] In general, it is also possible to construct monotone, subadditive, continuous from below set functions that are not integrable. \\ Consider  the set function given in  Example \ref{esempio1}.ii).
For each singleton \(\{n\}\) we have
\(\nu(\{n\})=\min \{1,\dfrac{1}{n}\}=
\dfrac{1}{n}.\)
Hence, for every $A \subset \mathbb{N}$ with $\card A = \aleph_0$, then 
\[
\nu(A) < \sum_{n \in A}\nu(\{n\}) \quad \mbox{and} \quad  (B_w)\int_{T} 1_A d\nu =
\sum_{n \in A}\frac{1}{n}=+\infty.
\]
\end{itemize}
\end{Example}

From now on, $\mathcal{M}_{cs} (\mathcal{E})$ is the set of all set functions $\nu:\mathcal{E}\rightarrow [0,\infty)$, with $\nu(\emptyset)=0$,
which are $\mathcal{E}$-integrable, continuous from below and subadditive.\\

In \cite{CG}, the following result was  given that shows the inequalities of
H\"{o}lder and Minkowski.
\begin{Theorem} \mbox{\rm (\cite[Theorems 3.8 and  3.9]{CG})}
\label{BUMI}
Let $\nu  \in \mathcal{M}_{cs} (\mathcal{E})$ and $u, v:T\rightarrow \mathbb{R}$ be measurable functions.
\begin{itemize}
\item[\ref{BUMI}.a)] If  $|u|^{p}, |v|^{q}, |uv| \in B_w(\nu,T)$, 
then
\begin{equation*}
\Vert u v\Vert_{1}\leq \Vert u\Vert_{p} \cdot \Vert v\Vert_{q}.
\quad  \text{ \rm (H\"{o}lder Inequality)}
\end{equation*}
\item[\ref{BUMI}.b)]
Suppose that $|u|^{p}, |v|^{p}, |u+v|^{p}, |u|\cdot|u+v|^{p-1}, |v|\cdot|u+v|^{p-1}\in B_w(\nu, T)$. Then
\begin{equation*}
  \Vert u+v\Vert_{p}\leq \Vert u\Vert_{p}+ \Vert v\Vert_{p}.
\quad \text{ \rm (Minkowski  Inequality)}
\end{equation*}
\end{itemize}
\end{Theorem}
 As in the classic case, the inequalities of H\"{o}lder and Minkowski are very important in the definition of the norm
of the spaces of integrable functions.
The following result is a consequence of Theorem \ref{Tix}.
\begin{Proposition}\label{modulo}
Let $u: T\rightarrow \mathbb{R}$ be a real function such that $u\in B_w(\nu)$. Then $|u|\in B_w(\nu,T).$
\end{Proposition}

Let \, $\mathcal{L}^{1}_{B_w}(\nu, T)$ be 
a linear subspace of 
$B_w(\nu, T) \cap \mathscr{F}(T, \mathbb{R})$.
As usually, denote by $L^{1}_{B_w}(\nu, T)$ the quotient space of $\mathcal{L}^{1}_{B_w}(\nu, T)$ with respect to the usual equivalence
relation "$\sim$":
\begin{equation*}
\mbox{for every } u, v\in \mathcal{L}^{1}_{B_w}(\nu, T), u\sim v  \textrm{ iff }  u= v \textrm{ }\nu \textrm{-ae.}
\end{equation*}
\begin{Theorem}
Suppose $\nu\in \mathcal{M}_{cs} (\mathcal{E})$. Then the function $\Vert \cdot\Vert _{1}$ is a norm on the space $L^{1}_{B_w}(\nu, T).$
\end{Theorem}
\begin{proof}
The proof is analogous to the classic one, using the properties of the Birkhoff weak integral of \cite{cgi,CG}.
\end{proof}
\begin{Remark}
\rm
 $B_w(\nu) \cap \mathscr{F}(T, \mathbb{R})$ is a subspace of $\mathcal{L}^{1}_{B_w}(\nu, T)$. In fact 
if $u \in B_w(\nu) \cap \mathscr{F}(T, \mathbb{R})$, then by construction and Proposition \ref{modulo}, for every $E \in \mathcal{E}$,  $u, \vert u \vert \in B_w(\nu, E)$. So if $u, v \in B_w(\nu)$, then $\alpha u + \beta v  \in B_w(\nu)$ by \cite[Theorema 4.3 and 4.5]{cgi} for every $\alpha, \beta \in \mathbb{R}$.
 Again by Proposition \ref{modulo}, $\vert \alpha u + \beta v \vert \in B_w(\nu)$.
\end{Remark}
Finally, if   $p\in (1, \infty)$, we can define analogously  the space $\mathcal{L}^{p}_{B_w}(\nu, T)$.\\

As said before, the Minkowski  and H\"{o}lder inequalities  fail when $0 < p < 1$.
We give now some examples in the non additive case.

\begin{Example}\label{esempio3} \rm 
Let $T = \mathbb{N}$ and  $\mathcal{E}$ be its power set. 
Take \(p\in(0,1)\) (e.g. \(p=\frac{1}{2}\)) and consider
\(
u=\mathbf{1}_{\{1\}},\, v=\mathbf{1}_{\{2\}}
\). The supports of \(u\) and \(v\) are disjoint.
\begin{itemize}
\item[\ref{esempio3}.i)] 
Let $\nu_1$ as in Examples \ref{esempio1}.i) and \ref{esempio2}.i)
and
\begin{eqnarray*}
&&
\|u\|_{p,\nu_1}^p=\nu_1(\{1\})=2^{-1}\quad
\|v\|_{p,\nu_1}^p=\nu_1(\{2\})=2^{-2},
\\ &&
\|u+v\|_{p,\nu_1}^p =\|u\|_{p,\nu_1}^p+\|v\|_{p,\nu_1}^p
=\frac{1}{2}+\frac{1}{4}=\frac{3}{4}.
\end{eqnarray*}
For  \(p=\dfrac{1}{2}\)  we have
\[
\big(\frac{3}{4}\big)^{2}= \frac{9}{16}
>
\big(\frac{1}{2}\big)^{2}+\big(\frac{1}{4}\big)^{2}=\frac{1}{4}+\frac{1}{16},
\]
so \(\|u+v\|_{p,\nu_1}\not\leq \|u\|_{p,\nu_1}+\|v\|_{p,\nu_1}\). \\
\item[\ref{esempio3}.ii)]
Let $\nu_2$ as in Examples 
\ref{esempio1}.ii) and 
\ref{esempio2}.ii).
Now
\begin{eqnarray*}
&& \|u\|_{p,\nu_2}^p=\nu_2(\{1\})=1,\qquad
\|v\|_{p,\nu_2}^p=\nu_2(\{2\})=\frac{1}{2},
\qquad 
\|u+v\|_{p,\nu_2}^p=\frac{3}{2}.
\end{eqnarray*}
For \(p=\dfrac{1}{2}\) this yields
\(\big(\dfrac{3}{2}\big)^{2} >
1^{2}+\big(\dfrac{1}{2}\big)^{2},\) \,
so again 
\[\|u+v\|_{p,\nu_2}\not\leq \|u\|_{p,\nu_2}+\|v\|_{p,\nu_2}\] and then the  Minkowski inequality fails.
\end{itemize}
\end{Example}

\begin{Remark}
\rm
The failure of the Minkowski inequality  is a consequence of the fact that 
 $p<1$ makes the map $t\mapsto t^p$ concave.
 What we are able to prove is that
 Minkowski's inequality reverses its direction under suitable assumptions when \(0<p<1\), 
 and therefore the usual subadditivity of the \(L^{p}\)-quasinorm fails and is replaced by a superadditivity inequality.
 \end{Remark}
 Regarding concave maps we can also observe that
 \begin{Lemma}\label{young-r}
 {\rm \cite[Corollary 2.1]{BB1}}
 Let \(0<p<1\) be fixed and set
\(q:=\frac{p}{p-1}<0.\)
Then, 
\[\,bc \;\ge\; \frac{b^{p}}{p}+\frac{c^{q}}{q}, \qquad \mbox{ for all }\,  b>0,\ c>0. \]
\end{Lemma}

We introduce now the main results of this paper: the 
 reverse H\"{o}lder's and Minkowski's inequalities for $0<p< 1$ and
 other inequalities for Birkhoff weak integrable scalar functions.
\begin{Theorem}\label{T1}
  Let $\nu \in \mathcal{M}_{cs} (\mathcal{E})$  and  $u, v \in \mathscr{F}(T, \mathbb{R})$.
  Let $p\in (0, 1)$ and $q$ is its conjugate. If 
\begin{description}
\item[\rm \ref{T1}.a)] \, $|uv|, |u|^{p}, |v|^{q}
\in B_w(\nu, T)$ and $\int_{T}|v|^{q}d\nu > 0$,  then
\begin{equation*}
\phantom{a} \qquad \Vert uv\Vert_{1}\geq \Vert u\Vert_{p} \cdot \Vert v\Vert_{q}.
\qquad  \text{ \rm (Reverse H\"{o}lder Inequality)}
\end{equation*}
\item[\rm \ref{T1}.b)] \,  $(|u|+|v|)^{p}, |u|^{p}$, $|v|^{p}, |u|(|u|+|v|)^{p-1}, |v|(|u|+|v|)^{p-1}$
$\in B_w(\nu, T)$,\, then
\begin{equation*}
\phantom{aa} \phantom{a} \qquad  \Vert\, |u|+|v| \, \Vert_{p}\geq  \Vert u \Vert_{p}+ \Vert v \Vert_{p}.
\qquad \text{ \rm (Reverse Minkowski  Inequality)}
\end{equation*}
\end{description}
\end{Theorem}

\begin{proof}
\begin{itemize}
\item[(\ref{T1}.a)] If $\int_{S}|u|^{p}d\nu =0$, then by \cite[Theorem 3.7]{CG}
  it results $uv=0$ $\nu-a.e.$ Therefore, the inequality of \ref{T1}.a) is true.\\
Suppose $\int_{S}|u|^{p}d\nu >0$.
Let 
   \begin{eqnarray*}
   b=  |u| \cdot \Big(\int_{T}|u|^{p}d\nu \Big)^{-\frac{1}{p}}, 
   \qquad 
c= |v| \cdot \Big(\int_{T}|v|^{q}d\nu \Big)^{-\frac{1}{q}}.
\end{eqnarray*}
Since for every $  b, c\in(0, \infty),$
 it is 
 $bc\geq \dfrac{b^{p}}{p}+\dfrac{c^{q}}{q}$, 
  by Lemma \ref{young-r},  then, in our setting, we have:
\begin{eqnarray*}
  \frac{|uv|}{(\int_{T}|u|^{p}d\nu)^{\frac{1}{p}}(\int_{T}|v|^{q}d\nu)^{\frac{1}{q}}}
  \geq
\frac{|u|^{p}}{p(\int_{T}|u|^{p}d\nu)}+ \frac{|v|^{q}}{q(\int_{T}|v|^{q}d\nu)}.
\end{eqnarray*}
  According to \cite[Theorems 5.5 and 6.1]{cgi} we have
\begin{eqnarray*}
 \frac{\int_{T}|uv|d\nu}{(\int_{T}|u|^{p}d\nu)^{\frac{1}{p}}(\int_{T}|v|^{q}d\nu)^{\frac{1}{q}}}
  &\geq&   \frac{\int_{T}|u|^{p}d\nu}{p \Big(\int_{T}|u|^{p}d\nu \Big)}+
\frac{\int_{T}|v|^{q}d\nu}{q \Big(\int_{T}|v|^{q}d\nu \Big)}=1
\end{eqnarray*}
and the assertion holds.
\item[(\ref{T1}.b)]
 From (\ref{T1}.a), it follows that:
\begin{eqnarray} \label{**}
  \int_{T}(|u|+|v|)^{p}d\nu &=&
\int_{T}(|u|+|v|)^{p-1}(|u|+|v|)d\nu\geq \\
&\geq& \Big(\int_{T}(|u|+|v|)^{q(p-1)}d\nu \Big)^{\frac{1}{q}}
\Big(\int_{T}|u|^{p}d\nu \Big)^{\frac{1}{p}}+ \nonumber\\
  &+& \Big(\int_{T}(|u|+|v|)^{q(p-1)}d\nu \Big)^{\frac{1}{q}}
\Big( \int_{T}|v|^{p}d\nu \Big)^{\frac{1}{p}}= \nonumber \\
&=& \nonumber
\Big( \int_{S}(|u|+|v|)^{p}d\nu  \Big)^{\frac{1}{q}}(\|u\|_{p}+ \|v\|_{p}).
\end{eqnarray}
Now we divide \, \eqref{**} \, by
$\Big(\displaystyle{\int_{T}}(|u|+|v|)^{p}d\nu \Big)^{\frac{1}{q}}$ and
 the Reverse Minkowski Inequality is obtained.
\end{itemize}
\end{proof}
Moreover
\begin{Corollary}\label{L1}
Let \(\nu \in \mathcal{M}_{cs}(\mathcal{E})\) and let \(u, v \in \mathscr{F}(T, \mathbb{R})\)
 with \(v(t)\neq 0\) for every \(t\in T\).
Let \(p\in (0,\infty)\setminus\{1\}\) and let \(q\) be its conjugate index.
Assume that 
\( |u|,   |v|, \frac{|u|^{p}}{|v|^{p/q}} \in B_w(\nu,T)\).
Then the form of the H\"{o}lder inequality depends on the value of \(p\):  
when \(p>1\) we have
\begin{eqnarray}\label{corL1-a}
   \Big(\int_{T} |u|\, d\nu \Big)^{p}
   \leq 
   \Big(\int_{T} \frac{|u|^{p}}{|v|^{p/q}} \, d\nu \Big)
   \cdot
   \Big(\int_{T} |v|\, d\nu \Big)^{p/q},
\end{eqnarray}
whereas for \(p\in(0,1)\) the inequality is reversed:
\begin{eqnarray}\label{corL1-b}
   \Big(\int_{T} |u|\, d\nu \Big)^{p}
   \geq 
   \Big(\int_{T} \frac{|u|^{p}}{|v|^{p/q}} \, d\nu \Big)
   \cdot
   \Big(\int_{T} |v|\, d\nu \Big)^{p/q}.
\end{eqnarray}

\end{Corollary}

\begin{proof}
 Assume $p>1$.  Applying Hölder's inequality, \ref{BUMI}.a), for conjugate indices  \(p,q\) to the functions
\( \dfrac{|u|}{|v|^{1/q}}\) and \( |v|^{1/q}.\)
Since, by hypotheses, 
\(\dfrac{|u|^{p}}{|v|^{p/q}}, \,\, |v|\in  B_w(\nu,T), \)
  we obtain
\[ \int_{T}|u|\,d\nu \le
\Big(\int_{T}\frac{|u|^{p}}{|v|^{p/q}}\,d\nu\Big)^{1/p}
\Big(\int_{T}|v|\,d\nu\Big)^{1/q}, \]
and raising both sides to the power \(p\) yields 
\[ \Big(\int_{T}|u|\,d\nu\Big)^{p} \le
\Big(\int_{T}\frac{|u|^{p}}{|v|^{p/q}}\,d\nu\Big)
\Big(\int_{T}|v|\,d\nu\Big)^{p/q}. \]
Applying  \ref{T1}.a) to the same functions, when $0 < p < 1$, we obtain the reverse inequality.
\end{proof}

\begin{Remark}\label{L2}
\rm
 Observe that, when $p> 1$,
\(|u|\in B_w(\nu,T)\) in Corollary \ref{L1}, is  a consequence of the two others  integrability assumptions via the inequality we prove.\\
 The hypothesis plays a role in the case $0 < p < 1$. \\ Moreover,  
by Corollary \ref{L1}, we obtain
  the following inequalities:
\begin{eqnarray*}
&&  \int_{T}\frac{|u|^{p}}{|v|^{p-1}}d\nu \geq \dfrac{\Big(\int_{T}|u|d\nu \Big)^{p}}{\Big(\int_{T}|v|d\nu \Big)^{p-1}}, \qquad \mbox{for every }  p>1;
\\
&&  \int_{T}\frac{|u|^{p}}{|v|^{p-1}}d\nu \leq
 \frac{\Big(\int_{T}|u|d\nu \Big)^{p}}{\Big(\int_{T}|v|d\nu \Big)^{p-1}},
 \qquad \mbox{for every }  p \in (0,1).
\end{eqnarray*}

In fact for  $p>1$,\, taking  $ \dfrac{|u|}{|v|^{1/q}}, \, |v|^{1/q}$, in  \eqref{corL1-a}, we have:
\[\int_T |u|\,d\nu 
\le  \Big(\int_T \frac{|u|^p}{|v|^{p-1}} \, d\nu\Big)^{1/p} \Big(\int_T |v| \, d\nu\Big)^{1/q}
\]
and then
\[\Big(\int_T |u|\,d\nu\Big)^p \le \Big(\int_T \frac{|u|^p}{|v|^{p-1}} \, d\nu\Big) \Big(\int_T |v| \, d\nu\Big)^{p-1},\]
While, 
for $0<p<1$,  applying \eqref{corL1-b}, we obtain
\[ \int_T |u|\,d\nu 
\ge  \Big(\int_T \frac{|u|^p}{|v|^{p-1}} \, d\nu\Big)^{1/p} \Big(\int_T |v| \, d\nu\Big)^{1/q}
\]
and so
\[ \Big(\int_T |u|\,d\nu\Big)^p \ge \Big(\int_T \frac{|u|^p}{|v|^{p-1}} \, d\nu\Big) \Big(\int_T |v| \, d\nu\Big)^{p-1}, \]
\end{Remark}
Finally
\begin{Theorem}\label{Tf}
Let $\nu \in \mathcal{M}_{cs} (\mathcal{E})$ and conjugate indices  $p, q\in (1, \infty)$. 
Suppose that  $u, v:T\rightarrow (0, \infty)$ are measurable  functions
and that there exist $\alpha, \beta \in(0, \infty)$ such that:
\begin{itemize}
\item[\rm \ref{Tf}.a)]
$\alpha\leq \dfrac{u(t)}{v(t)}\leq \beta$,\,   for every $ t\in T.$
   If $u, v, u^{\frac{1}{p}}v^{\frac{1}{q}}\in B_w(\nu, T)$,\,  then
\begin{eqnarray*}
  \Big( \int_{T}ud\nu \Big)^{\frac{1}{p}} \cdot
\Big( \int_{T}vd\nu \Big)^{\frac{1}{q}}\leq
 \Big(\frac{\beta}{\alpha} \Big)^{\frac{1}{pq}}\cdot
  \int_{T}u^{\frac{1}{p}}v^{\frac{1}{q}}d\nu.
\end{eqnarray*}
\item[\rm \ref{Tf}.b)]
$\alpha\leq \dfrac{u^{p}(t)}{v^{q}(t)}\leq \beta$,\,  for every  $t\in T.$
   If $uv$, $u^{p}$, $v^{q}\in B_w(\nu, T)$,\,  then
\begin{eqnarray*}
  \Big(\int_{T}u^{p}d\nu \Big)^{\frac{1}{p}} \cdot
\Big(\int_{T}v^{q}d\nu \Big)^{\frac{1}{q}}\leq
\Big(\frac{\beta}{\alpha} \Big)^{\frac{1}{pq}}\cdot
  \int_{T}u  v\, d\nu.
\end{eqnarray*}
\end{itemize}
\end{Theorem}
\begin{proof}
\begin{itemize}
\item[\ref{Tf}.a)] For every 
$ t\in T$, it is  $\dfrac{u(t)}{v(t)}\leq \beta$, therefore
\[  v^{\frac{1}{q}}(t)\geq \beta^{-\frac{1}{q}}u^{\frac{1}{q}}(t).\]
Then we have, for every $t \in T$, 
\begin{equation*}
 u^{\frac{1}{p}}(t)\, v^{\frac{1}{q}}(t)\geq \beta^{-\frac{1}{q}}u^{\frac{1}{p}}(t)\, u^{\frac{1}{q}}(t)=
 \beta^{-\frac{1}{q}}\, u(t).
\end{equation*}
By 
\cite[Theorems 5.5 and 6.1]{cgi}
 it follows that
\begin{equation}\label{formula1}
 \Big(\int_{T}u^{\frac{1}{p}}\, v^{\frac{1}{q}}d\nu \Big)^{\frac{1}{p}}\geq \beta^{-\frac{1}{pq}}
 \Big(\int_{T}ud\nu \Big)^{\frac{1}{p}}.
\end{equation}
Since $\dfrac{u(t)}{v(t)}\geq \alpha$, for every  $ t\in T,$ we have
$  u^{\frac{1}{p}}(t)\geq \alpha^{\frac{1}{p}}v^{\frac{1}{p}}(t)$
and
\begin{equation*}
 u^{\frac{1}{p}}(t)\, v^{\frac{1}{q}}(t)\geq \alpha^{\frac{1}{p}}v^{\frac{1}{p}}(t)\, v^{\frac{1}{q}}(t)=
 \alpha^{\frac{1}{p}}\,  v(t),\qquad   \mbox{for every } t\in T.
\end{equation*}
 By \cite[Theorems 5.5 and 6.1]{cgi}, it results
\begin{equation}\label{formula2}
 \Big(\int_{T}u^{\frac{1}{p}}\, v^{\frac{1}{q}}d\nu \Big)^{\frac{1}{q}}\geq \alpha^{\frac{1}{pq}}
 \Big(\int_{T}vd\nu \Big)^{\frac{1}{q}}.
\end{equation}
According to (\ref{formula1}) and (\ref{formula2}), the inequality of \ref{Tf}.a) follows.
\item[\ref{Tf}.b)]
 It holds if we consider $u^{p}$ and $v^{q}$ instead of $u$ and $v$ in \ref{Tf}.a).
\end{itemize}
\end{proof}

\section{Applications}\label{sec4}
In this section we quote some  applications and some  future fields of research.

\subsection{Vector valued case}
 We can extend our result to vector functions $u: T \to X$ where $X$ is a Banach space. The definition of $B_w -\nu$ integrability is the same as in Definition \ref{def-int}, where we consider the $\Vert  \cdot \Vert _X$ instead of $\vert \cdot \vert$. In this case we will use the symbol $B^X_w (\nu, T)$.  Some results are already obtained for what concernes integrability and convergence results, \cite{cgi,CGI15,CG}. 
\\
Let $X$ be a Banach space, $p \in (0,+\infty[$ and $u:T \to X$ be a function satisfying the condition $\Vert u \Vert_X^p \in B_w(\nu, T)$.
Denote 
\[ \Vert u \Vert_{X,p} = \Big( \int_T \Vert u \Vert_X^p d \nu \Big)^{1/p}. \]
The following result can be obtained
\begin{Theorem}\label{vector}
Consider $\nu \in \mathcal{M}_{cs} (\mathcal{E})$  and  $u, v: T \to X$ measurable functions, 
  $p\in (0, 1)$ and its conjugate  $q$. Suppose
\begin{itemize}
\item[\rm \ref{vector}.a)] \, $\Vert u \Vert_X \cdot \Vert v \Vert_X, 
\Vert u \Vert_X^{p}, \Vert v \Vert_X^{q}
\in B_w(\nu, T)$ and $\int_{T}\Vert v \Vert_X^{q}d\nu > 0$,  then
\begin{equation*}
\phantom{a} \qquad \int_T \Vert u \Vert_X \cdot \Vert v \Vert_X d \nu \geq \Vert u\Vert_{X,p} \cdot \Vert v\Vert_{X,q}.
\qquad  \text{ \rm (Reverse H\"{o}lder Inequality)}
\end{equation*}
\item[\rm \ref{vector}.b)] \,  
$(\Vert u \Vert_X +\Vert v\Vert_X)^{p}, 
\Vert u \Vert_X^{p} \cdot \Vert v \Vert_X ^{p}, 
\Vert u \Vert_X (\Vert u \Vert_X +\Vert v\Vert_X)^{p-1}, 
\Vert v\Vert_X(\Vert u \Vert_X +\Vert v \Vert_X)^{p-1}
\in B_w(\nu, T)$,\, then
\begin{equation*}
\phantom{aa} \phantom{a} \qquad 
\big \vert\, \Vert u \Vert_X +\Vert v\Vert_X \, \big\vert_{p}\geq  \Vert u \Vert_{X,p}+ \Vert v \Vert_{X,p}.
\qquad \text{ \rm (Reverse Minkowski  Inequality)}
\end{equation*}
\end{itemize}
\end{Theorem}
\begin{Remark}
\rm
Similar results to Corollary \ref{L1}, Remark \ref{L2} and Theorem \ref{Tf} can be obtained.
\end{Remark}
\subsection{Applications in  Interval Analysis}
 An important field with many applications is Interval Analysis. In 1966, Moore \cite{Mo} used for the first time elements of Interval Analysis in numerical analysis and computer science. Interval-valued functions have many applications in uncertainty theory, signal and image processing or in edge detection algorithms (e.g. \cite{CCGIS20,mesiar1,ZZZ,W}).
\subsubsection{Interval-valued functions and scalar set functions ($F,\nu$)}

In \cite{Medjm,CGI15} the authors defined the Birkhoff weak (simple) integral of multifunctions and presented some of its properties making use of the Hausdoff distance and of the  R{\aa}dstr\"{o}m embedding.
Integral inequalities of interval-valued functions were obtained for example in \cite{CGIS} with respect to different
types of integrability. Integral inequalities are important tools in computing deviations or measuring actions.\\

Let $\big( ck(\mathbb{R}_0^+), \oplus, \cdot, d_H, \preceq \big)$ be the  complete (not linear)  metric space (see for example \cite{flores}) consisting of 
  the family of all non-empty bounded and closed intervals  of 
  $\mathbb{R}_0^+$  with the   Minkowski addition  $\oplus$,  the multiplication by scalars
   and the Hausdorff distance $d_H$, where $d_H(A;B) :=\sup_{x\in A}\,\inf_{y\in B}\,|x-y|$. 
   The hyperspace $ck([0,1])$  is used, for example, in decision theory.
\begin{itemize}
    \item 
    By convention  $\{0\}=[0,0]$
and  $\Vert A\Vert :=\sup\{|x|\colon x\in{A}\}$.
\item In  particular, if $A=[a^-,a^+] \in ck(\mathbb{R}_0^+)$ then $\Vert A\Vert =a^+$,  $d_H([0,a^+], [0,b^+]) = |b^+ -a^+|$  and 
 $d_H([a^-,a^+], [b^-,b^+]) = \max \{|a^- -b^-|,|a^+ - b^+|\}$. Moreover \(d_H(A\oplus B,\;C \oplus D) \leq d_H(A,C)+d_H(B,D).\)
 \item
Finally  $\preceq$ is the weak interval order, namely 
$[a^-,a^+] \preceq [b^-,b^+]$ if and only if $a^- \leq b^-$ and $a^+ \leq b^+$.    
\item 
 In order to work in this setting, and in particular to study the H\"{o}lder inequality,  we consider also the following 
 operation: 
\[ [a^-,a^+] \otimes [b^-, b^+] = [a^- \cdot b^-, a^+ \cdot b^+]
=[a^- b^-, a^+ b^+].\]
 \end{itemize}
 Following \cite{CCGIS}, let $F: T \to  ck(\mathbb{R}_0^+)$ be the interval-valued function 
 defined by $F(t) = [u^-(t), u^+(t)]$  with $u^-(t) \leq u^+(t)$
 for every $t \in T$. The functions $u^{\pm}$ are the so-called endpoints functions of the interval $F(t)$. In this setting we can also consider the multifunctions determined by integrable functions, as in \cite{CDMS2019}.

\begin{Definition}\label{def31}
\rm
 $F:T\to ck(\mathbb{R}_0^+)$ is said to be 
 Birkhoff weakly
 integrable (on $T$) with respect to $\nu$, if there is an interval $I \in ck(\mathbb{R}_0^+)$
 such that for every $\varepsilon>0$, there exist
$P_{\varepsilon}\in \mathcal{C}$ and $n_{\varepsilon}\in \mathbb{N}$ 
such that for every $  P \in \mathcal{C}$,\,  $P=(B_{n})_{n\in \mathbb{N}}$,\, $P\geqslant P_{\varepsilon}$
and every $t_{n}\in B_{n}, n\in \mathbb{N}$:
\begin{eqnarray*}\label{defi}
d_H (\oplus_{k=1}^{n}F(t_{k})\nu(B_{k}), I) <\varepsilon, \quad \mbox{ for every } n\geq n_{\varepsilon}.
\end{eqnarray*}
The interval $I$ is called
the Birkhoff weak integral of $F$ on $T$ with respect to $\nu$ and is denoted by  $I:=\int_{T} F d\nu$.
Obviously, if it exists, is unique.
\end{Definition}
 Analogously to  the scalar case, we denote by  $\mathbf{B}_w(\nu,T)$  the family of all interval-valued functions that are integrable in the Birkhoff weak sense on $T$ with respect to $\nu$.
 
\begin{Remark}\label{rem-somma}\rm  
Given $  P=(B_{k})_{k\in \mathbb{N}} \in \mathcal{C}$, and    $t_{k}\in B_{k}$ we said that  $ \{(B_{k}, t_k)\}_{k\in \mathbb{N}}$ is a tagged partition. For   every $n \in \mathbb{N}$, we have 
\begin{eqnarray*}
\sigma_{F,n} (P) &:=& \oplus_{k=1}^{n} F(t_k) \nu(B_k) 
=  \oplus_{k=1}^{n} [u^-t (t_k), u^+(t_k)] \nu(B_k)= \\
&=& \{ \sum_{k=1}^{n} y_k, \,\, y_k \in
  [u^- (t_k) , u^+(t_k) ]\nu(B_k) \}.
\end{eqnarray*}
Since the set $\sigma_{F,n} (P) \in ck(\mathbb{R}_0^+)$, so it is an interval 
$[u_{n}^{-^{(P)}}, u_{n}^{+^{(P)}}]$.\\
\end{Remark}

	\begin{Theorem}\label{p2}	
	 $F =[u^-, u^+] \in \mathbf{B}_w(\nu,T)$ if and only if $u^{\pm}
	\in B_w(\nu,T)$   and
	\begin{eqnarray}\label{form-1}
 \int_T F d\nu = \left[ \int_T u^-(t) d\nu, \,  \int_T u^+(t)  d\nu\right]. 
	\end{eqnarray}
\end{Theorem}
\begin{proof} Let  $F =[u^-, u^+] \in \mathbf{B}_w(\nu,T)$. So there exists $I=[a,b] \in ck(\mathbb{R}_0^+)$ such that  for every $\varepsilon>0$, there exist $n_{\varepsilon} \in \mathbb{N}$ and  a countable partition
$P_{\varepsilon}$ of $T$, so that for every tagged partition
$P=\{(A_{n}, t_n)\}_{n\in \mathbb{N}}$ of $T$ with $P\geqslant P_{\varepsilon}$,
by Remark \ref{rem-somma}, we have:
\[d_H ([u_{n}^{-^{(P)}}, u_{n}^{+^{(P)}}], [a,b]) :=
 \max \big\{
\vert  u_{n}^{-^{(P)}} - a \vert\, , \, 
\vert  u_{n}^{+^{(P)}} - b \vert \big\}
 <\varepsilon, \mbox{ for every } n\geq n_{\varepsilon}.\]
So, in particular, 
for every tagged partition
$P=\{(A_{n}, t_n)\}_{n\in \mathbb{N}}$ of $T$ with $P\geqslant P_{\varepsilon}$ and for every  $n\geq n_{\varepsilon}$, 
 we have:
\[\max \Big\{ \big\vert \sum_{k=1}^{n} u^- (t_k) \nu  (A_k) - a \big\vert ,\, 
\big\vert \sum_{k=1}^{n} u^+(t_k) \nu  (A_k) - b \big\vert \, \Big\}
 \leq \varepsilon,\]
 and then  $u^{\pm} \in B_w(\nu,T)$. 
The equality \eqref{form-1} follows from the definition of the integral.\\
  For the converse implication, for every $\varepsilon>0$, let 
  $n_{\varepsilon}>0$ and  let
  $ P_{\varepsilon}$  be a countable partition that verify the definition of  Birkhoff weak integrability  on $T$  for both $u^{\pm}$.
 Then, for every partition
$P:=\{ B_n, \, n \in \mathbb{N} \} \geqslant P_{\varepsilon}$, for every $t_n \in B_n$ we have, for every $n \geq n_{\varepsilon}$, 
\begin{eqnarray*}
\max \Big\{ 
\big\vert \sum_{k=0}^{n} u^-(t_{k}) \nu  (B_{k})- \int_T u^- d\nu \big \vert, \,
\big\vert \sum_{k=0}^{n} u^+(t_{k}) \nu  (B_{k})- \int_T u^+ d\nu \big \vert
\Big\}
<\varepsilon.
\end{eqnarray*}
Since $u^{\pm}(t) \in F(t)$,  for all $t \in T$, this implies
\[d_H \left([u_{n}^{-^{(P)}}, u_{n}^{+^{(P)}}],
 \left[ \int_T u^- d\nu, \, \int_T u^+  d\nu\right] \right)
  \leq \varepsilon
  \mbox{ for every } n \geq n_{\varepsilon}\]
and then the assertion follows.
\end{proof}
Moreover the  Birkhoff weak integrability
is hereditary on subsets  $A\in\mathcal{E}$. In fact
\begin{Proposition}\label{intersection}
 $F \in \mathbf{B}_w(\nu,A)$ for every $A \in \mathcal{E}$
 if and only if $F\,  1_{A} \in \mathbf{B}_w(\nu,T)$
and 
 \[\int_{A} F\, d\nu =\int_T F \,1 _{A}\, d\nu .\]
\end{Proposition}

\begin{proof}
Assume that $F \in \mathbf{B}_w(\nu,T)$ and denote by $[a,b]$ its integral. \\
Then  for every $\varepsilon>0$,  there exists a countable partition $P_{\varepsilon}$ of $T$, such that, for every finer countable partition
 $P:=\{B_{n}\}_{n\in \mathbb{N}}$  and for every $t_n \in B_n$ we have 
\[d_H \left( \sigma_{F,n} (P), [a,b] \right) \leq \varepsilon, \qquad 
 \mbox{for every } n \geq n_{\varepsilon}.\]
Let $A\in \mathcal{E}$ and 
 $P_0 \geqslant P_{\varepsilon} \wedge  \{A, T\setminus A\}$.
Let  $P_A \subset P_0$ the  partition 
of  $A$.\\
 Let $\Pi_A \geqslant P_A$ be a 
 partition of $A$ and extend it with a common partition of $T\setminus A$  in such a way the new partition is finer
than $P_{\varepsilon}$. \\
It is possible to prove that
$\sigma_{F,n}(\Pi_A)$
is  Cauchy  in $ck(\mathbb{R}_0^+)$, and so the first claim follows by the completeness of the space.
The equality follows from  Proposition \ref{p2}.
\end{proof}
Analogously to \cite{CCGIS} and using Theorem  \ref{p2} and Proposition \ref{intersection}, 
 it is possible to prove some properties of the Birkhoff weak integral of an interval-valued function: 
\begin{itemize}
\item 
The Birkhoff weak integral is additive with respect to the Minkowski addition. In fact, given $F := [u^-,u^+]$ and $G := [v^-,v^+]$ in  $\mathbf{B}_w(\nu,T)$, we can write
\begin{eqnarray*}
\int_T \big( F \oplus G) d\nu &=&
\int_T [ u^-(t)+v^-(t), u^+(t)+v^+(t)] d\nu = 
\\ &=& \Big[ \int_T (u^-(t) + v^-(t)) d\nu, \int_T (u^+(t) + v^+(t)) d\nu \Big] =
\\ &=&
\Big[ \int_T u^-(t)  d\nu, \int_T u^+(t)  d\nu \Big] \oplus
\Big[ \int_T v^-(t)  d\nu, \int_T v^+(t)  d\nu \Big] = \\
&=& \int_T F  d\nu \oplus \int_T G d\nu.
\end{eqnarray*}
\item 
Using the same notation as before, if 
$p > 0$ we have
\[
\big(F\oplus G\big)^p \;=\;
\big([u^- +v^-,\; u^+ + v^+]\big)^p
\;=\;
\big[(u^-+v^-)^p,\;(u^+ + v^+)^p\big].
\]
In fact, 
on \(\mathbb{R}_0^{+}\) the map \(x\mapsto x^p\) is  increasing.
Since \(F\oplus G=[u^- + v^-,\;u^+ +v^+]\), the image of \(F\oplus G\) under the map \(x\mapsto x^p\) is the interval
\[
\big\{x^p:\;x\in[u^- + v^-,\;u^+ +v^+]\big\}=
\big[(u^- +v^-)^p,\;(u^+ +v^+)^p\big].
\].
\item
If $F^p \in \mathbf{B}_w(\nu,T)$, then 
\[
\Big\{\int_T  F^p d\nu \Big\}^{1/p}
:=\Big[ \Big(\int_T (u^-)^p\,d\nu\Big)^{1/p},\;\Big(\int_T (u^+)^p\,d\nu\Big)^{1/p}\Big]
= \Big[ \Vert u^-\Vert _{p},\;\Vert u^+\Vert _{p}\Big].
\]
We denote by 
\[ \Vert F\Vert _p:= d_H \Big( \Big\{\int_T  F^p d\nu \Big\}^{1/p}, \{0\} \Big)
= \Vert u^+\Vert _{p}.\]
\begin{eqnarray*}
\Big\{\int_T \big( F \oplus G \big)^p d\nu \Big\}^{1/p} &=& 
\Big\{ \int_T [ (u^- + v^-)^p, (u^+ + v^+)^p] d\nu \Big\}^{1/p}= 
\\ &=&
 \Big[ \Big\{ \int_T (u^- + v^-)^p d\nu\Big\}^{1/p} , \Big\{ \int_T (u^+ + v^+)^p] d\nu \Big\}^{1/p}\Big].  
\end{eqnarray*}
\end{itemize}
So, if we assume for $u^{\pm}, v^{\pm}$ the hypotheses of Theorem \ref{BUMI}.b) or Theorem \ref{T1}.b), we are able to obtain a Minkowski or a reverse Minkowski inequality, depending on the value of $p > 0$.
\begin{Theorem}
 Let $\nu \in \mathcal{M}_{cs} (\mathcal{E})$  and let 
 $F = [u^-,u^+]$ and $G = [v^-,v^+]$ with 
 $u^{\pm}, v^{\pm} \in 
 \mathscr{F}(T, \mathbb{R}_0^+)$.
  Let $p >0$. If 
 $ ( u^{\pm})^{p}, \, 
  ( v^{\pm} )^{p}, \,
 (u^- + u^+)^{p}, (v^- + v^+)^{p}, \, 
 u^{\pm}(u^- + u^+ )^{p-1}, \,
 v^{\pm}(v^- +v^+)^{p-1}
\in B_w(\nu, T)$, 
then the Minkowski, for $p \geq 1$ (or the reverse Minkowski, for $0 < p < 1$)  inequality holds.
\end{Theorem}
\begin{proof} 
By hypotheses $F \oplus G \in \mathbf{B}_w(\nu,T)$ and 
\begin{eqnarray*}
\Vert F \oplus G \Vert_p &=& 
\Big\{\int_T \big( F \oplus G)^p d\nu \Big\}^{1/p} = 
\Big\{ \int_T [ (u^- + v^-)^p, (u^+ + v^+)^p] d\nu \Big\}^{1/p}= 
\\ &=& 
\Big[ \Big\{ \int_T (u^- + v^-)^p d\nu\Big\}^{1/p} , 
\Big\{ \int_T (u^+ + v^+)^p] d\nu \Big\}^{1/p}\Big] = \\
&=& \Big[ \Vert u^- + v^- \Vert_p\, ,\, \vert u^+ + v^+ \Vert_p \Big]  
\end{eqnarray*}
We split now the proof in two cases, depending on the values of $p$.
\begin{description}
\item[Case $p \geq 1$)]
By the  Minkowski inequality given in Theorem \ref{BUMI}.b), applied to $u^{\pm}, v^{\pm}$,  we have
\[
\Vert u^- + v^-\Vert _{p} \le \Vert u^-\Vert _{p}+\Vert v^-\Vert _{p},
\qquad
\Vert u^+ + v^+\Vert _{p} \le \Vert u^+\Vert _{p}+\Vert v^+\Vert _{p}.
\]
and so, in $ck(\mathbb{R}_0^+)$ with the weak interval order we have:
\begin{eqnarray*}
 \Big[ \Vert u^- + v^-\Vert _{p},\;\Vert u^+ + v^+\Vert _{p}\Big]
&\leqslant &
\Big[ \Vert u^-\Vert _{p}+\Vert v^-\Vert _{p},\;\Vert u^+\Vert _{p}+\Vert v^+\Vert _{p}\Big]=\\
&=& 
\Big[ \Vert u^-\Vert _{p},\;\Vert u^+\Vert _{p}\Big]
\oplus
\Big[\Vert v^-\Vert _{p},\;\Vert v^+\Vert _{p}\Big];  \\
\Vert F \oplus G \Vert_p &=&  \Vert u^+ + v^+\Vert _{p} \leq 
\Vert u^+\Vert _{p}+\Vert v^+\Vert _{p} = \\
&=& \Vert F \Vert_p + \Vert G \Vert_p.
\end{eqnarray*}
\item[Case $0 < p < 1$)]
We proceed in the same way, this time using the reverse Minkowski inequality given in Theorem \ref{T1}.b), applied to $u^{\pm},v^{\pm}$.
\end{description}
\end{proof}
We are ready now to examine the H\"{o}lder inequality for $p > 1$.
 \begin{Theorem}
Let $\nu  \in \mathcal{M}_{cs} (\mathcal{E})$  and let 
 $F = [u^-,u^+]$ and $G = [v^-,v^+]$ with  $u^{\pm},\, v^{\pm} \in 
 \mathscr{F}(T, \mathbb{R}_0^+)$
 be  such that 
  $(u^{\pm})^{p},\,  (v^{\pm})^{q}, (u^- u^+), \, (v^- v^+) \in B_w(\nu,T)$.
  Then, if $p \in \mathbb{R}^+ \setminus \{1 \}$ and $q$ is the coniugate index, we have
  \begin{eqnarray*} 
 && \Vert F \otimes G \Vert_1 \preceq \Vert F \Vert_p \cdot \Vert G \Vert_q, \quad \mbox{ when } p > 1;\\
 && \Vert F \otimes G \Vert_1 \succeq \Vert F \Vert_p \cdot \Vert G \Vert_q, \quad \mbox{ when } 0 < p < 1 \mbox{\,\, and  }  \int_T (u^+)^q d\nu > 0,  \int_T (v^+)^q d\nu > 0.
  \end{eqnarray*}
\end{Theorem}
\begin{proof} 
Let $F = [u^-,u^+]$ and $G = [v^-,v^+]$.
By the definition of the integral, the definition of $\otimes$, the weak interval relation $\preceq$, Theorems \ref{p2} and  \ref{BUMI}.a) we have:
\begin{eqnarray*}
\int_T F \otimes G d\nu &=& \int_T \big[u^- v^-, u^+ v^+\big] d\nu = 
\Big[ \int_T u^- v^- d\nu, \int_T u^+ v^+ d\nu \Big] =\\
&=& \Big[ \Vert u^-  v^- \Vert_1, \Vert u^+ v^+ \Vert_1 \Big] \leqslant \Big[ \Vert u^-\Vert_p \cdot \Vert v^- \Vert_q, \Vert u^+\Vert_p \cdot \Vert v^+ \Vert_q \Big]\\
\Vert F \otimes G \Vert_1 &\leq& \Vert u^+ \Vert_p \cdot \Vert v^+ \Vert_q =
\Vert F \Vert_p \cdot \Vert G \Vert_q.
\end{eqnarray*}
For the second inequality we apply  Theorem \ref{T1}.a) instead of Theorem \ref{BUMI}.a) to the components of $F$ and $G$.
\end{proof}

\subsubsection{Scalar functions and interval-valued set functions ($u,M$)}
\phantom{a} \\
$M$ is said to be an {\it interval submeasure} if 
$M: \mathcal{E} \to ck(\mathbb{R}_0^+)$ is defined by 
\(M(A) = [\nu^-(A), \nu^+(A)]\),
where $\nu^{\pm}: \mathcal{E} \to \mathbb{R}_0^+$
are two monotone and subadditive set functions  with $\nu^-(A) \leq \nu^+(A)$ for every $A \in \mathcal{E}$.
An interval submeasure $M$ satisfies:
\begin{itemize}
\item $M (\emptyset )=\{0\};$
\item $M(A) \preceq M(B)$ for every $A, B \in \mathcal{E}$ with $A \subseteq B$ \,\, (monotonicity);
\item $M(A \cup B) \preceq M(A) \oplus M(B)$ for every  $A, B \in \mathcal{E}$ with $A \cap B = \emptyset$\,\, (subadditivity).
\end{itemize}
An example of interval-valued submeasure  was  given in \cite{de,sh}, where a mathematical theory of evidence was proposed using  the non-additive measures called  Belief and Plausibility;
  in such a way for every set $A$ un interval was associated:   
  \[[\nu^-(A), \nu^+(A)] :=[Bel(A),Pl(A)].\]
    For results in this subject see for example \cite{Pap1,Pap}.\\
  
 $M$ is an  additive multimeasure  if $M(A \cup B) = M(A) \oplus M(B)$ for every disjoint sets $A, B \in \mathcal{E}$.
By \cite[Remark 3.6]{Pap}, $M(A) = [\nu^-(A),\, \nu^+(A)]$ is a multisubmeasure with respect to ``$\preceq$'' if and only if $\nu^{\pm}$
 are monotone and subadditive.\\
 
Now the Birkhoff weak integrability of a scalar function with respect to a interval-valued multisubmeasure follows easily:
 \begin{Definition}
\rm
 $u:T\to \mathbb{R}$ is said to be 
is Birkhoff weakly
 integrable (on $T$) with respect to $M: \mathcal{E} \to ck(\mathbb{R}_0^+)$, if there exists  an interval $I \in ck(\mathbb{R}_0^+)$
 such that for every $\varepsilon>0$, we can find
$P_{\varepsilon}\in \mathcal{C}$ and $n_{\varepsilon}\in \mathbb{N}$ 
such that for every $  P \in \mathcal{C}$ with $P=(B_{n})_{n\in \mathbb{N}}$ and $P\geqslant P_{\varepsilon}$,
and every $t_{n}\in B_{n}, n\in \mathbb{N}$ we have 
\begin{eqnarray*}\label{defi}
d_H (\oplus_{k=1}^{n}u(t_{k})M(B_{k}), I) <\varepsilon, \quad \mbox{ for every } n\geq n_{\varepsilon}.
\end{eqnarray*}
The interval $I:=\int_{T} u\, dM$ and is called
the Birkhoff weak integral of $F$ on $T$ with respect to $M$.
Obviously, if it exists, is unique.\\
 Let $\mathbf{B}_w(M,T)$ be  the family of all scalar functions that are integrable in the Birkhoff weak sense on $T$ with respect to $M$. 
\end{Definition}
\begin{Remark}
\rm It is easy to see that if $u \in \mathbf{B}_w(M,T)$, then 
\[ \int_T u\, dM = \Big[ \int_T u\,  d\nu^-, \int_T u \, d\nu^+ \Big] \]
and
\begin{eqnarray*}
&& d_H \Big(\oplus_{k=1}^{n}u(t_{k})M(B_{k}), 
\Big[ \int_T u \, d\nu^-, \int_T u\,  d\nu^+ \Big] \Big) = \\ &&=
\max \Big\{ \big\vert \sum_{k=1}^{n}u(t_{k})\nu^-(B_{k}) -
\int_T u \, d\nu^- \big\vert, \,
\big\vert \sum_{k=1}^{n}u(t_{k})\nu^+(B_{k}) -
\int_T u\,  d\nu^+ \big\vert \Big\}.
\end{eqnarray*}
So the linearity follows from the scalar-scalar ($u,\nu$) case, in fact:
\begin{eqnarray*}
\int_T (u+v) \, dM &=& \Big[\int_T (u+v) d\nu^-, \int_T (u+v) d\nu^+ \Big] =\\ &=&
\Big[\int_T u\,  d\nu^- + \int_T v\, d\nu^-, 
\int_T u\,  d\nu^+ + \int_T v\,  d\nu^+ \Big] = \\
&=& \int_T u\,  dM \oplus \int_T v \, dM.
\end{eqnarray*}
\end{Remark} 
At this point the Minkowski inequality or its reverse are obtained, as in the interval-valued - scalar case ($F,\nu$), thanks to the scalar-scalar ($u,\nu$) case. 
For the H\"{o}lder inequality we can observe that
if the hypotheses of the Theorem \ref{BUMI}.a) or Theorem \ref{T1}.a) are verified, then
\begin{eqnarray*}
    \int_T u(t) v(t) \, dM &=& \Big[ \int_T u(t) v(t) \, d\nu^-, 
    \int_T u(t) v(t)\, d\nu^+\Big] \quad \mbox{and}\\
    \Big \Vert \int_T u(t) v(t) \, dM  \Big\Vert &=&
    \Big \vert \int_T u(t) v(t) \, d\nu^+ \Big\vert,
\end{eqnarray*}
i.e.,
\begin{eqnarray*}
    \Vert uv \Vert_{1,M} &=& \Vert uv \Vert_{1,\nu^+},
\end{eqnarray*}
and so H\"{o}lder and reverse H\"{o}lder inequalities immediately follow  from the scalar-scalar ($u,\nu$) case, when $\nu=\nu^+$.
\begin{Remark}\label{re-ir1}
\rm
Interval-valued functions provide a natural modeling framework in image processing
when uncertainty, quantization effects, or measurement errors cannot be neglected.
In fact, 
representing images as interval-valued functions allows one to explicitly encode this uncertainty at each spatial location, rather than treating it as an  additive perturbation.
For example, in applications such as fractal image coding for compression and edge detection, representing images as interval-valued objects allows one to explicitly account for
round--off errors and discretization inaccuracies arising in the conversion of
 analog signals into digital data.
More generally, when pixel values cannot be assigned with full precision,
an interval-valued approach becomes preferable.
This situation is common in images acquired through physical measurement processes,such as biomedical imaging (e.g.\ CT or MR imaging), where sensor sensitivity and data conversion introduce intrinsic uncertainty.
In this setting, interval-valued functions offer a unified formalism to encode
noise and imprecision at the pixel level, enabling robust image analysis and processing methods that preserve relevant structural information while explicitly handling uncertainty.
For results in this subject
see for example \cite{latorre, mesiar1,CCGIS20,CS2024}.
\end{Remark}


\vspace{6pt}

\subsection{A  sparsity-promoting reconstruction model }\label{sec4.3}
We can model a gray--scale image as a measurable function
\(
 f : T \subset \mathbb{R}^2 \to \mathbb{R},
\) such that $f^p \in B_w(T,\nu)$.
In practical acquisition systems, the observed image is often corrupted by noise,
which can be written as
\(|g(x)-f(x)| \le |\eta(x)|\) \,  for a.e.  \(x \in T\),
where \(\eta^p \in B_w(T,\nu)\) represents a noise and  $g : T \to \mathbb{R}$
denotes the observed image, i.e.\ the data provided by the acquisition
process.\\
This condition is natural in many acquisition models, including:
$g = f + \eta$, \,
$g = f + \psi(f,\eta)$ \  with  $|\psi(f,\eta)| \le |\eta|$.
A  sparsity-promoting reconstruction model
promotes compact representation since  only a few coefficients carry significant information; preserves  edges, textures, and fine details are maintained while suppressing noise.

While standard models are known in the cases of functions belonging to Lebesgue spaces $L^p$,  it can be useful to consider a 
  sparsity-promoting reconstruction model in the quasi-Banach
regime \(0<p<1\)\, as the variational problem in $B_w(T, \nu)$ given by

\[
\min_{f^p \in B_w(T, \nu)}
\left\{
\mathcal{D}(f,g)
\;+\;
\lambda \int_\Omega |f(x)|^p \, d\mu(x)
\right\},
\]

where \(\mu\) is a submeasure modeling spatial interaction or uncertainty;
\(\mathcal{D}(f,g)\) is a data fidelity term that  measures the discrepancy between the candidate reconstruction \(f\) and the
observed data \(g\). In a non-additive setting, a natural choice could be a
Birkhoff weak-type functional of the form
\[
\mathcal{D}(f,g)
:=
\int_\Omega |f(x)-g(x)|^r \, d\mu(x),
\qquad r>0;
\] 
\(\lambda>0\) is a regularization parameter and it balances fidelity and sparsity.
The regularization functional
\[
\mathcal{R}(f) := \int_\Omega |f(x)|^p \, d\mu(x), \qquad 0<p<1,
\]
is nonconvex and strongly promotes sparsity, in the sense that it penalizes
small nonzero coefficients more heavily than convex \(B_w(T, \nu)\)-type penalties.
In other words, the $0<p<1$ regularization naturally produces images with few nonzero components, which is useful in problems such as compression, denoising, or reconstruction from incomplete data.
\\
In general, $r$ and $p$ are  independent.
For example $r=2$ corresponds to a 
least-squares type fidelity (classical in additive noise). 
For sparsity-promoting reconstruction models in the standard setting $L^p$  see for example \cite{KSB2024,PST2025}.
\section{Conclusion}
We have proved several inequalities and reverse forms of the classical H\"{o}lder and Minkowski
inequalities for Birkhoff weakly integrable functions when the underlying set function with
respect to which we integrate is non-additive. Moreover, these results have also been
extended to  vector-valued integrands and to the 
multivalued setting, where the integrands, or the set functions involved, 
are allowed to take interval-valued outputs.\\
H\"{o}lder's and Minkowski's inequalities remains an open problem when both the function $F$ and the set functions $M$ are interval-valued. 
In this latter setting, results concerning integrability, convergence, and properties of
the integral have been established, while the corresponding inequalities are currently the
subject of ongoing investigation. Also applications in Image  Reconstruction, in the spirit of Remark \ref{re-ir1} and of Section \ref{sec4.3}, will be exploited in a future study.
\\


\vspace{6pt} 



 {\small 
{\bf Funding}
This research was partly funded,
for the last two authors,
 by the Unione Europea - Next Generation EU, Missione 4 Componente C2 - CUP Master: J53D2300390 0006, CUP: J53D23003920 006 - Research project of MUR (Italian Ministry of University and Research) PRIN 2022  “Nonlinear differential problems with applications to real phenomena” (Grant Number: 2022ZXZTN2);
 PRIN 2022 PNRR: “RETINA: REmote sensing daTa INversion with multivariate functional modelling for essential climAte variables characterization” funded by the European Union under the Italian National Recovery and Resilience Plan (NRRP) of NextGenerationEU, under the MUR (Project Code: P20229SH29, CUP: J53D23015950001)
 and INdAM - GNAMPA Project   2024 Dyna.M.I.Ch.E "Dynamical Methods: Inverse problems, Chaos and Evolution",  CUP Id: E53C23001670001.
\\

{\bf Acknowledgments} This research has been accomplished by the last two authors  within the
 UMI Group TAA - “Approximation Theory and Applications”,
 the G.N.AM.P.A. group of INDAM, the University of Perugia and the ``Centro  di Ricerca Interdipartimentale Lamberto Cesari'' of  the University of Perugia.\\

The authors declare no conflicts of interest.\\
}

\Addresses


\begin{thebibliography}{999}
\small
\bibitem{BB1}
Benaissa, B., 
On the reverse Minkowski's integral inequality
\textit{Kragujevac Journal of Mathematics}, 
46 (3) (2022), 407--416.

\bibitem{BB2}
Benaissa, B., 
A generalization of reverse H\"{o}lder's inequality via the diamond-$\alpha$
integral on time scale, 
\textit{Hacet. J. Math. Stat.}, 51(2) (2022), 383--389.

\bibitem{Frey2022}
 Bernicot, F., Frey,  D.,
Characterizations of weak reverse Hölder inequalities on metric measure spaces,  
\textit{Mathematische Zeitschrift}, 301:2269–2290, (2022).  

\bibitem{BS}
Boccuto, A., Sambucini, A. R., 
A note on comparison between Birkhoff and Mc Shane integrals for multifunctions,
\textit{ Real Analysis Exchange}, 37 (2) (2012), 3-15.

\bibitem{BCS}
Boccuto, A., Candeloro, D., Sambucini, A.R.,
 Henstock multivalued integrability in Banach lattices with respect to pointwise non atomic measures,
\textit{ Atti Accad. Naz. Lincei Rend. Lincei Mat. Appl.}, 26 (4), (2015), 363-383.

\bibitem{BS3}
Boccuto, A., Sambucini, A.R., Abstract integration with respect to measures and applications to modular convergence in vector lattice setting, 
\textit{Results in Mathematics}, (2023) 78:4.

\bibitem{C}
Candeloro, D., 
Riemann-Stieltjes integration in Riesz spaces,  
\textit{Rend. Mat. Appl.}, (7) 16(4), 563–585 (1996).

\bibitem{cdms2020}
Candeloro, D., Di Piazza, L., Musia\l, K. Sambucini, A. R., 
Multi-integrals of finite variation, 
 {\it Boll. Unione Mat. Ital.},  13, 459–468 (2020). 

\bibitem{CS}
Candeloro, D., Sambucini A.R., 
Comparison between some norm and order gauge integrals in Banach lattices,
\textit{Pan American Math. J.}, 
25 (3) (2015), 1-16.

\bibitem{CCGS}
Candeloro, D., Croitoru, A., Gavrilu\c{t}, A., Sambucini, A.R., 
Atomicity related to non-additive integrability, 
\textit{Rend. Circ. Mat. Palermo, II. Ser.}, 65 (2016), 435–449.

\bibitem{Medjm}
Candeloro, D., Croitoru, A., Gavrilu\c{t}, A., Sambucini, A.R., 
 An extension of the Birkhoff integrability for multifunctions, 
\textit{Mediterr. J. Math.}, 
13(5), 2551–2575 (2016). 

\bibitem{cms}
Candeloro, D., Mesiar, R.,  Sambucini, A. R.,
A special class of fuzzy measures: Choquet integral and applications,
\textit{Fuzzy Sets and Systems}, 355 (2019) 83 -- 99.

\bibitem{CDMS2019}
Candeloro, D.,  Di Piazza, L., Musia\l, K., Sambucini, A.R.,
Multifunctions determined by integrable functions,
{\it International Journal of Approximate Reasoning},
112, (2019),  140–148.

\bibitem{CCGIS}
Candeloro, D., Croitoru, A. Gavrilu\c{t}, A., Iosif, A., Sambucini, A.R.,
Properties of the Riemann-Lebesgue integrability in the non-additive case,
\textit{Rend. Circ. Mat. Palermo, II. Ser.}, 69 (2020), 577–589.

\bibitem{CCGIS20}
Costarelli, D., Croitoru, A., Gavrilu\c{t}, A., Iosif, A., Sambucini, A. R.,
The Riemann-Lebesgue integral of interval-valued multifunctions,
\textit{Mathematics}, 
8 (2020), 1--17.

\bibitem{CS2024}
Costarelli, D., Sambucini, A.R., 
A comparison among a fuzzy algorithm for image rescaling
with other methods of digital image processing,
{\it Constructive Mathematical Analysis}, 
7 (2), (2024), 45-68.

\bibitem{CC}
Chen, G., Chen, Z.,
A functional generalization of the reverse H\"{o}lder integral inequality on time scale,
\textit{Math. Comput. Model.}, 54 (2011), 2939--2941.

\bibitem{CCS}
Chichon K., Chichon M., Satco B., 
Differential inclusions and multivalued integrals,
\textit{Discuss. Math. Differ. Incl. Control Optim.}, 
33 (2013), 171–191.

\bibitem{CGI15}
Croitoru, A., Gavrilu\c{t}, A., Iosif, A., 
The Birkhoff weak integrability of multifunctions,
\textit{International Journal of Pure Mathematics}, 2 (2015), 47--54.

\bibitem{cgi}
Croitoru, A., Gavrilu\c{t}, A., Iosif, A., 
The Birkhoff weak integral of functions relative to a set function in Banach spaces setting integrability,
\textit{WSEAS Transaction on Mathematics},  16 (2017), 375--383.

\bibitem{CG}
Croitoru, A., Gavrilu\c{t}, A., 
Convergence results in Birkhoff weak integrability,
\textit{Bollettino dell’Unione Matematica Italiana},
13 (2020), 477--485.

\bibitem{CGIS}
Croitoru, A. Gavrilu\c{t}, A., Iosif, A., Sambucini, A.R., Inequalities in Riemann-Lebegue integrability,
\textit{Mathematics},
12(1) (2024), 49.

\bibitem{cis2025}
Croitoru,  A., Iosif, A., Sambucini, A.R.,
Generalized Decomposition Integral of Real Functions with Respect to Fuzzy Measures,
\textit{Bol. Soc. Paran. Mat.},
43 (2), (2025), 1-13.

\bibitem{de}
 Dempster, A. P., Upper and lower probabilities induced by a multivalued mapping, 
{\it Ann. Math.Statist.}, { 1967}, { 38}, 325--339.

\bibitem{Feller}
Feller, W., \textit{An Introduction to Probability Theory and its Applications}, vol. 2,
John Wiley and Sons, Inc. New York, 1966.

\bibitem{flores}
Roman-Flores, H., Chalco-Cano, Y., Lodwick, W. A.,
Some integral inequalities for interval-valued functions,
\textit{Comp. Appl. Math.}, 
37 (2018), 1306--1318. 

\bibitem{GIC}
Gavrilu\c{t}, A., Iosif A., Croitoru, A., 
The Gould integral in Banach lattices, 
\textit{Positivity},
19 (1) (2015), 65-82.

\bibitem{GCPS}
Giurgescu  M., Chi\c{t}escu, I., Paraschiv, T., \c{S}tefan, C.,
\textit{Diagnostication method of openness using the non-linear integrals}, U.P.B. Sci. Bull., Series A, 82 (1) (2020), ISSN 1223-7027.


\bibitem{mesiar1} 
  Jurio, A.,  Paternain, D.,  Lopez-Molina,  C.,   Bustince, H.,   Mesiar, R. and  Beliakov, G.,
   A construction method of
interval-valued Fuzzy Sets for image processing, 2011 IEEE Symposium on Advances in Type-2 Fuzzy Logic Systems
(T2FUZZ), (2011), 16–22.

\bibitem{KHB}
Kalita, H., Hazarika, B., Becerra, T. P., 
On AP-Henstock–Kurzweil Integrals and Non-Atomic Radon Measure,
\textit{Mathematics},
11, (2023), 1-17.

\bibitem{Kawabe}
Kawabe J., Yamada N.,
The completeness and separability of function spaces in non additive measure theory,
\textit{Fuzzy Sets and Systems}, 466 (2023), 108409.

\bibitem{KSB2024}
 Kumar, N.,  Sonkar, M.,  Bhatnagar, G.,
Efficient image restoration via non-convex total variation regularization and ADMM optimization, {\it Applied Mathematical Modelling},
 132, (2024),  428-453, ISSN 0307-904X,
Doi: 10.1016/j.apm.2024.04.055.

\bibitem{latorre} 
   La Torre, D., Mendivil, F., 
   Minkowski-additive multimeasures, monotonicity and self-similarity,
{\it Image Analysis and Stereology}, { (2011)}, {\it 30} (3),  135 -- 142,
Doi:10.5566/ias.v30.p135-142.

\bibitem{Milne}
Milne, E. A., Note on Rosseland’s integral for the stellar absorption coefficient, \textit{Monthly Notices R Astron.Soc.}, 85 (1925), 979--984.

\bibitem{Mo}
Moore, R. E., \textit{Interval Analysis}, Prentice Hall, Englewood Cliffs, NJ, USA, 1966.

\bibitem{MS2024} 
Marraffa, V., Sambucini, A.R.,
Vitaly theorems for varying measures
\textit{Symmetry}, 
 2024, 16, 972. 

\bibitem{radko}
Mesiar, R., Li, J., Pap, E.,
The Choquet Integral as Lebesgue integral and related inequalities,
\textit{Kybernetika},  46 (6) (2010), 1098--1107.

\bibitem{Pap1}
Pap, E., Pseudo-additive measures and their applications, in: E. Pap (Ed.), Handbook of
Measure Theory, II, Elsevier, { 2002}, 1403--1465.

\bibitem{Pap}
Pap, E., Iosif, A., Gavrilu\c{t}, A., 
Integrability of an Interval-valued Multifunction with respect to an Interval-valued Set Multifunction, 
\textit{ Iranian Journal of Fuzzy Systems},
15 (3)(2018), 47–63.


\bibitem{PS}
Precupanu, A., Satco, B., 
The Aumann-Gould Integral, 
\textit{Mediterr. J. Math.}, 5 (2008), 429–441.

\bibitem{PST2025}
Prater-Bennette, A., Shen, L., Tripp,    Wei, J.,
 Algorithms for structured sparsity promoting functions regularized image restoration model,  
 {\it Sampl. Theory Signal Process. Data Anal.},  23, 10 (2025). Doi: 10.1007/s43670-025-00102-7

\bibitem{sh}
 Shafer, G., {\it A Mathematical Theory of Evidence}, Princeton University Press, Princeton, { 1976}.

\bibitem{SC}
Stamate C., Croitoru A., 
\textit{Non-linear integrals, properties and relationships} in:
Recent Advances in Telecommunications, Signals and Systems (Proceedings of NOLASC 13), WSEAS Press, 2013, 118--123.

\bibitem{SC25}
Stamate, C.,  Croitoru, A., 
Integral Inequalities for Vector (Multi)functions,
\textit{Axioms}
2025, 14(12), 915 (2025).

\bibitem{torra}
Torra, V., Narukawa, Y., Sugeno, M. (eds), 
\textit{Non-additive measures: Theory and applications},
Studies in Fuzziness and Soft Computing 310 (2014), Springer, Switzerland.

\bibitem{W}
Wohlberg, B., De Jager, G., 
A review of the fractal image coding literature,
\textit{IEEE Trans. Image Process.},  8 (1999)
1716–1729.

\bibitem{Yang2023}
 Yang W.,
Certain New Reverse H\"{o}lder- and Minkowski-Type Inequalities for Modified Unified Generalized Fractional Integral Operators with Extended Unified Mittag–Leffler Functions,
\textit{Fractal and Fractional}, 7(8):613, (2023).

\bibitem{Yin}
Yin L., Qi F., 
Some Integral Inequalities on Time Scales, 
\textit{ Results. Math.}, 64, 371–381 (2013). 

\bibitem{ZZZ}
Zhou, Y., Zhang, C., Zhang, Z., 
An efficient fractal image coding algorithm using unified feature and DCT,
\textit{Chaos Solitons Fractals},
 39 (2009), 1823–1830.

\end{thebibliography}
\end{document}